\documentclass[12pt]{amsart}
\usepackage{amscd,amssymb}
\usepackage{amsthm,amsmath,amssymb}
\usepackage[matrix,arrow]{xy}

\newtheorem{theorem}[equation]{Theorem}

\newtheorem{proposition}[equation]{Proposition}
\newtheorem{lemma}[equation]{Lemma}
\newtheorem{corollary}[equation]{Corollary}
\newtheorem{conjecture}[equation]{Conjecture}

\theoremstyle{definition}
\newtheorem{example}[equation]{Example}
\newtheorem{definition}[equation]{Definition}

\theoremstyle{remark}
\newtheorem{remark}[equation]{Remark}

\makeatletter\@addtoreset{equation}{section} \makeatother

\author{Ivan Cheltsov}

\title{Points in projective spaces and applications}

\address{\begin{tabbing}
\hspace*{28 em}\=\kill
Steklov Institute of Mathematics \>School of Mathematics\\
8 Gubkin street, Moscow 117966   \>University of Edinburgh\\
Russia                           \> Edinburgh EH9 3JZ, UK\\
                                 \>\\
\texttt{cheltsov@yahoo.com}      \>\texttt{I.Cheltsov@ed.ac.uk}
\end{tabbing}}

\thanks{We assume that all varieties are projective, normal, and defined over $\mathbb{C}$.}%

\begin{document}

\begin{abstract}
We prove the factoriality of a nodal hypersurface in
$\mathbb{P}^{4}$ of degree $d$ that has at most $2(d-1)^{2}/3$
singular points, and factoriality of a double cover of
$\mathbb{P}^{3}$ branched over a nodal surface of degree $2r$
having less than $(2r-1)r$ singular points.
\end{abstract}

\maketitle

\section{Introduction.}
\label{section:intro}

Let $\Sigma$ be a finite subset in $\mathbb{P}^{n}$ and
$\xi\in\mathbb{N}$, where $n\geqslant 2$. The points of the set
$\Sigma$ impose independent linear conditions on homogeneous forms
of degree $\xi$ if and only if for every point $P\in\Sigma$ there
is a homogeneous form of degree $\xi$ that vanishes at
$\Sigma\setminus P$ and does not vanish at $P$, which is
equivalent to
$h^{1}(\mathcal{I}_{\Sigma}\otimes\mathcal{O}_{\mathbb{P}^{n}}(\xi))=0$,
where $\mathcal{I}_{\Sigma}$ is the ideal sheaf of $\Sigma$.

In this paper we prove the following result (see
Section~\ref{section:main}).

\begin{theorem}
\label{theorem:main} Suppose that at most $\lambda k$ points of
the set $\Sigma$ lie on a curve of degree~$k$, where
$\lambda\in\mathbb{N}$ and $\lambda\geqslant 2$. Then
$h^{1}(\mathcal{I}_{\Sigma}\otimes\mathcal{O}_{\mathbb{P}^{n}}(\xi))=0$
if one of the following conditions~holds:
\begin{itemize}
\item $\xi=\lfloor 3\lambda/2-3\rfloor$ and $|\Sigma|<\lambda\lceil\lambda/2\rceil$;%
\item $\xi=\lfloor 3\mu-3\rfloor$ and $|\Sigma|\leqslant\lambda\mu$, where $\mu\in\mathbb{Q}$ such that $\lfloor3\mu\rfloor-\mu-2\geqslant\lambda\geqslant\mu$;%
\item $\xi=\lfloor n\mu\rfloor$ and $|\Sigma|\leqslant\lambda\mu$, where $\mu\in\mathbb{Q}$ such that $(n-1)\mu\geqslant\lambda$.%
\end{itemize}
\end{theorem}

Let us consider applications of Theorem~\ref{theorem:main}.

\begin{definition}
\label{definition:factorial} An algebraic variety is called
factorial if its divisor class group is $\mathbb{Z}$.
\end{definition}

Let $\pi\colon X\to\mathbb{P}^{3}$ be a double cover branched over
a surface $S$ of degree $2r\geqslant 4$ such that the only
singularities of $S$ are isolated ordinary double points.  Then
$X$ is a hypersurface
$$
w^{2}=f_{2r}(x,y,z,t)\subset\mathbb{P}\big(1^{4},r\big)\cong\mathrm{Proj}\Big(\mathbb{C}[x,y,z,t,w]\Big),%
$$
where
$\mathrm{wt}(x)=\mathrm{wt}(y)=\mathrm{wt}(z)=\mathrm{wt}(t)=1$,
$\mathrm{wt}(w)=r$, and $f_{2r}$ is a homogeneous~po\-ly\-no\-mial
of degree $2r$ such that $f_{2r}=0$ defines the surface
$S\subset\mathbb{P}^{3}\cong\mathrm{Proj}(\mathbb{C}[x,y,z,t])$.

It follows from \cite{Ha70} and \cite{Di90} that the following
conditions are equivalent:
\begin{itemize}
\item the threefold $X$ is factorial;%
\item the singularities of the threefold $X$ are $\mathbb{Q}$-factorial;%
\item the equality $\mathrm{rk}\,H_{4}(X,\mathbb{Z})=1$ holds;%
\item the ring $\mathbb{C}[x,y,z,t,w]\slash I$ is a UFD, where $I=<w^{2}-f_{2r}(x,y,z,t)>$;%
\item the points of the set $\mathrm{Sing}(S)$ impose independent
linear conditions on homogeneous forms on $\mathbb{P}^{3}$ of
degree $3r-4$.%
\end{itemize}

In the case $r=3$, the threefold $X$ is known to be non-rational
if it is factorial (see \cite{ChPa04}), but the threefold $X$ is
rational if the surface $S$ is the Barth sextic (see \cite{Ba96}).

\begin{theorem}
\label{theorem:double-solid} Suppose that
$|\,\mathrm{Sing}(S)|<(2r-1)r$. Then $X$ is factorial.
\end{theorem}

\begin{proof}
The subset $\mathrm{Sing}(S)\subset\mathbb{P}^{3}$ is a
set-theoretic intersection of surfaces of degree $2r-1$, which
implies that $X$ is factorial by Theorem~\ref{theorem:main}.
\end{proof}

The claim of Theorem~\ref{theorem:double-solid} is conjectured in
\cite{Ch04t}, and it is proved in \cite{ChPa04} in the case $r=3$.

\begin{example}
\label{example:non-factorial} Suppose that the surface $S$ is
given by an equation
\begin{equation}
\label{equation:nodal-surface}
g^{2}_{r}(x,y,z,w)=g_{1}(x,y,z,w)g_{2r-1}(x,y,z,w)\subset\mathbb{P}^{3}\cong\mathrm{Proj}\Big(\mathbb{C}[x,y,z,w]\Big),%
\end{equation}
where $g_{i}$ is a general homogeneous polynomial of degree $i$.
Then $X$ is not factorial,~sin\-gu\-lar points of the surface $S$
are isolated ordinary double points, and
$|\,\mathrm{Sing}(S)|=(2r-1)r$.
\end{example}

We prove the following result\footnote{The claim of
Theorem~\ref{theorem:Hong-and-Park} is conjectured in
\cite{HoPa04}, and it is proved in \cite{HoPa04} the case $r=3$.}
in Section~\ref{section:extra}.

\begin{theorem}
\label{theorem:Hong-and-Park} Suppose that
$|\,\mathrm{Sing}(S)|\leqslant(2r-1)r+1$. Then the threefold $X$
is not~factorial if and only if the surface
$S\subset\mathbb{P}^{3}$ can be defined by the
equation~\ref{equation:nodal-surface}.
\end{theorem}

Let $V$ be a hypersurface in $\mathbb{P}^{4}$ of degree $d$ such
that $V$ has at most isolated ordinary double points. Then $V$ can
be given by the equation
$$
f_{n}(x,y,z,t,u)=0\subset\mathbb{P}^{4}\cong\mathrm{Proj}\Big(\mathbb{C}[x,y,z,t,u]\Big),%
$$
where $f_{n}$ is a homogeneous polynomial of degree $n$. It
follows from \cite{Ha70} and~\cite{Di90}~that~the~hy\-per\-surface
$V$ is factorial if and only is one of the following conditions
holds:
\begin{itemize}
\item the hypersurface $V$ has $\mathbb{Q}$-factorial singularities;%
\item the equality $\mathrm{rk}\,H_{4}(V,\mathbb{Z})=1$ holds;%
\item the ring $\mathbb{C}[x,y,z,t,u]/I$ is a UFD, where $I=<f_{n}(x,y,z,u)>$;%
\item the points of the set $\mathrm{Sing}(V)$ impose independent linear conditions on homogeneous forms on $\mathbb{P}^{4}$ of degree $2d-5$.%
\end{itemize}

In the case $d=4$, the hypersurface $V$ is not rational if it is
factorial (see \cite{Me03}), but a~ge\-ne\-ral determinantal
quartic threefold is rational and has isolated ordinary
double~points.

\begin{conjecture}
\label{conjecture:factoriality} Suppose that
$|\mathrm{Sing}(V)|<(d-1)^{2}$. Then $V$ is factorial.
\end{conjecture}

The claim of Conjecture~\ref{conjecture:factoriality} is proved in
\cite{Ch04e} and \cite{ChPa05d} in the case when $d\leqslant 7$.

\begin{example}
\label{example:non-factorial-hypersurface} Suppose that the
hypersurface $V$ is given by the equation
$$
xg(x,y,z,w,t)+yf(x,y,z,w,t)=0\subset\mathbb{P}^{4}\cong\mathrm{Proj}\Big(\mathbb{C}[x,y,z,w,t]\Big),%
$$
where $g$ and $f$ are general homogeneous polynomials of degree
$d-1$. Then $V$ is not fac\-to\-rial, singular points of $V$ are
isolated ordinary double points, and
$|\,\mathrm{Sing}(V)|=(d-1)^{2}$.
\end{example}

The factoriality of $V$ is proved in \cite{Ch04t} in the case when
$|\,\mathrm{Sing}(V)|\leqslant (d-1)^{2}/4$.

\begin{theorem}
\label{theorem:hypersurface} Suppose that
$|\,\mathrm{Sing}(V)|\leqslant 2(d-1)^{2}/3$. Then $V$ is
factorial.
\end{theorem}

\begin{proof}
The set $\mathrm{Sing}(V)$ is a set-theoretic intersection of
hypersurfaces of degree $d-1$, which implies the claim for
$d\geqslant 7$ by Theorem~\ref{theorem:main}. In the case
$d\leqslant 6$, the~factoriality of the hypersurface $V$ follows
from Theorem~2 in \cite{EiJ87} .
\end{proof}

Let $Y$ be a complete intersection of hypersurfaces $F$ and $G$ in
$\mathbb{P}^{5}$ of degree $m$ and $k$, respectively, such that
$m\geqslant k$ and $Y$ has at most isolated ordinary double
points.

\begin{example}
\label{example:complete-intersection-II} Let $F$ and $G$ be
general hypersurfaces that contain a two-dimensional linear
subspace in $\mathbb{P}^{5}$. Then $F$ and $G$ are smooth, the
threefold $Y$ has isolated ordinary double points, and
$|\mathrm{Sing}(Y)|=(m+k-2)^{2}-(m-1)(k-1)$, but $Y$ is not
factorial.
\end{example}

It follows from \cite{Di90} that the threefold $Y$ is factorial if
its singular points impose independent linear conditions on
homogeneous forms on $\mathbb{P}^{5}$ of degree $2m+k-6$.

\begin{theorem}
\label{theorem:complete-intersection-I}%
Suppose that $G$ is smooth, and $Y$ has at most
$(m+k-2)(2m+k-6)/5$ ordinary double points. Then the complete
intersection $Y$ is factorial for $m\geqslant 7$.
\end{theorem}

\begin{proof}
The set $\mathrm{Sing}(Y)$ is a set-theoretic intersection of
hypersurfaces of degree $m+k-2$, which concludes the proof by
Theorem~\ref{theorem:main}.
\end{proof}

Arguing as in the proof of
Theorem~\ref{theorem:complete-intersection-I}, we obtain the
following result.

\begin{theorem}
\label{theorem:complete-intersection-II}%
Suppose that $G$ is smooth, and $Y$ has at most
$(2m+k-3)(m+k-2)/3$ ordinary double points. Then the complete
intersection $Y$ is factorial for $m\geqslant k+6$.
\end{theorem}

Let $H$ be a smooth hypersurface in $\mathbb{P}^{4}$ of degree
$d$, and $\eta\colon U\to H$ be a double cover ramified in a
surface $R\subset H$ that is cut out by a hypersurface of degree
$2r\geqslant d$~that has isolated ordinary double points. Then $U$
is factorial if the points of $\mathrm{Sing}(R)$ impose
independent linear conditions on homogeneous forms of degree
$3r+d-5$ (see \cite{Di90}).

\begin{theorem}
\label{theorem:double-hypersurface}%
The threefold $U$ is factorial if
$|\mathrm{Sing}(R)|\leqslant(2r+d-2)r/2$ and $r\geqslant d+7$.
\end{theorem}

\begin{proof}
The set $\mathrm{Sing}(R)$ is a set-theoretic intersection of
hypersurfaces of degree $2r+d-2$, which implies the claim by
Theorem~\ref{theorem:main}.
\end{proof}

The author thanks I.\,Aliev, A.\,Corti, M.\,Gri\-nen\-ko,
V.\,Iskov\-s\-kikh, J.\,Park, Yu.\,Pro\-kho\-rov, V.\,Sho\-ku\-rov
and K.\,Shramov for useful and helpful conversations.

\section{Main result.}
\label{section:main}

Let $\Sigma$ be a finite subset in $\mathbb{P}^{n}$, where
$n\geqslant 2$. In this section we prove the following special
case of Theorem~\ref{theorem:main} leaving other cases to the
reader, because their proofs~are~similar.

\begin{proposition}
\label{proposition:3r-4} Suppose that at most $(2r-1)k$ points of
the set $\Sigma$ lie on a curve~of~degree~$k$, and
$|\Sigma|<(2r-1)r$, where $r\in\mathbb{N}$ and $r\geqslant 2$.
Then the points of the set $\Sigma$ impose independent linear
conditions on homogeneous forms of degree $3r-4$.
\end{proposition}

We may assume that $n\geqslant 3$ due to the following result,
which is Corollary~4.3 in \cite{DaGe88}.

\begin{theorem}
\label{theorem:Bese} Let $\pi\colon Y\to\mathbb{P}^2$ be a blow up
of points $P_{1},\ldots, P_{\delta}$, and $E_{i}$ be the
$\pi$-ex\-cep\-ti\-onal divisor such that $\pi(E_{i})=P_{i}$. Then
the linear system
$|\pi^{*}(\mathcal{O}_{\mathbb{P}^2}(\xi))-\sum_{i=1}^{\delta}E_{i}|$~does~not
have base points if at most $k(\xi+3-k)-2$ points of the set
$\{P_{1},\ldots, P_{\delta}\}$ lie on a curve of degree $k$ for
every natural number $k\leqslant (\xi+3)/2$, and the inequality
$$
\delta\leqslant
\mathrm{max}\left\{\Big\lfloor\frac{\xi+3}{2}\Big\rfloor\Big(\xi+3-\Big\lfloor\frac{\xi+3}{2}\Big\rfloor\Big)-1,
\Big\lfloor\frac{\xi+3}{2}\Big\rfloor^{2}\right\},%
$$
holds, where $\xi$ is a natural number such that $\xi\geqslant 3$.
\end{theorem}

Hence, to prove Proposition~\ref{proposition:3r-4}, we may assume
that $n=3$ due to the following result.

\begin{lemma}
\label{lemma:zero-dimensional} Let $\Pi$ be an $m$-dimensional
linear subspace in $\mathbb{P}^n$ such that $n>m\geqslant 2$,~and
$$
\psi\colon\mathbb{P}^n\dasharrow\Pi\cong\mathbb{P}^m
$$
be a projection from a general $(n-m-1)$-dimensional linear
subspace $\Omega\subset\mathbb{P}^{n}$ such that there is a subset
$\Lambda\subset\Sigma$ such that $|\Lambda|\geqslant\lambda k+1$,
but the set $\psi(\Lambda)$ is contained in an irre\-du\-ci\-ble
curve of degree $k$, and $\mathcal{M}$ be the linear system of
hypersurfaces in $\mathbb{P}^n$~of~degree~$k$~that contain
$\Lambda$. Then the base locus of $\mathcal{M}$ is
zero-di\-men\-si\-onal, and either $m=2$,~or~$k>\lambda$.
\end{lemma}

\begin{proof} Suppose that the base locus of $\mathcal{M}$ contains an
irreducible curve $Z$. Let $\Xi$ be a subset of the set $\Lambda$
consisting of points that are contained in the curve $Z$. Then we
may assume that $\psi(Z)$ does not contain any point of the set
$\psi(\Lambda\setminus\Xi)$, and $\psi\vert_{Z}$ is a birational
morphism, because $\psi$ is a general projection. Thus, we have
$\mathrm{deg}(\psi(Z))=\mathrm{deg}(Z)$.

Let $C$ be an irreducible curve in $\Pi$ of degree $k$ that
contains $\psi(\Lambda)$, and $W$ be the cone in $\mathbb{P}^{n}$
over the curve $C$ and with vertex $\Omega$. Then $W\in
\mathcal{M}$, which implies that $W$ contains the curve $Z$. Thus,
we have $\psi(Z)=C$, which implies that $\Xi=\Lambda$ and
$\mathrm{deg}(Z)=k$, but the curve $Z$ contains at most $\lambda
k$ points of the set $\Sigma$. Hence, the the base locus of the
linear system $\mathcal{M}$ is zero-dimensional.

Suppose that $m>2$ and $k\leqslant \lambda$. Let us show that the
latter assumption leads to a contradiction. We may assume that
$m=3$ and $n=4$, because we may consider $\psi$ as a composition
of $n-m$ projections from points. Thus, the projection
$\psi\colon\mathbb{P}^{4}\dasharrow\mathbb{P}^{3}$ is a projection
from the point $\Omega\in\mathbb{P}^{4}$.

Let $\mathcal{Y}$ be the set of all irreducible reduced surfaces
in $\mathbb{P}^{4}$ of degree $k$ that contains all points of the
set $\Lambda$, and $\Upsilon$ be a subset of $\mathbb{P}^{4}$
consisting of points that are contained in every surface of
$\mathcal{Y}$. Then $\Lambda\subseteq\Upsilon$, but the previous
arguments imply that $\Upsilon$ is a finite set.

Let $\mathcal{S}$ be the set of all surfaces in $\mathbb{P}^{3}$
of degree $k$ such that $S\in\mathcal{S}$ if and only if there is
a surface $Y\in\mathcal{Y}$ such that $\psi(Y)=S$ and
$\psi\vert_{Y}$ is a birational morphism. Then $\mathcal{S}$ is
not empty, because the projection $\psi$ is general enough and the
construction of the set $\mathcal{Y}$ does not depend on the
choice of the projection $\psi$. Let $\Psi$ be a subset of
$\mathbb{P}^{3}$ consisting of points that are contained in every
surface of the set $\mathcal{S}$. Then
$\psi(\Lambda)\subseteq\psi(\Upsilon)\subseteq\Psi$ by
construction.

The generality of $\Omega$ implies that $\psi(\Upsilon)=\Psi$.
Indeed, for every point $O\in\Pi\setminus\Psi$ and any general
surface $Y\in\mathcal{Y}$, we may assume that the line passing
through $O$ and $\Omega$ does not intersect $Y$, but the
restriction $\psi\vert_{Y}$ is a birational morphism.

Thus, the set $\Psi$ is a set-theoretic intersection of surfaces
in $\Pi$ of degree $k$, which implies that at most $\delta k$
points in $\Psi$ lie on a curve in $\Pi$ of degree $\delta$.
Hence, at most $k^{2}$ points of the set $\Psi$ lie on a curve in
$\Pi$ of degree $k$, but $\psi(\Lambda)$ contains at least
$\lambda k+1$ points contained in an irreducible curve in $\Pi$ of
degree $k$, which is a contradiction.
\end{proof}

Thus, we have a finite subset $\Sigma\subset\mathbb{P}^{3}$ such
that $|\Sigma|<(2r-1)r$, and at most $(2r-1)k$ points of $\Sigma$
lie on a curve of degree $k$, where $r\in\mathbb{N}$ and
$r\geqslant 2$. Fix an integer $\epsilon$ such that
$$
\big|\Sigma\big|<\big(2r-1\big)\big(r-\epsilon\big),
$$
and $\epsilon\geqslant 0$. We prove the following result, which
implies Proposition~\ref{proposition:3r-4}.

\begin{proposition}
\label{proposition:3r-4-epsilon} The points of the set $\Sigma$
impose independent linear conditions on homogeneous forms of
degree $3r-4-\epsilon$.
\end{proposition}

Fix an arbitrary point $P$ of the set $\Sigma$. To prove
Proposition~\ref{proposition:3r-4-epsilon} it is enough to
construct a surface\footnote{For simplicity we consider
homogeneous forms on $\mathbb{P}^{3}$ as surfaces in
$\mathbb{P}^{3}$.} in $\mathbb{P}^{3}$ of degree $3r-4-\epsilon$
that contains  $\Sigma\setminus P$ and does not contain $P$.

We may assume that $r\geqslant 3$ and $\epsilon\leqslant r-3$,
because the claim of Proposition~\ref{proposition:3r-4-epsilon}
follows from Theorem~2 in \cite{EiJ87} and
Theorem~\ref{theorem:Bese} in the case when $r\leqslant 3$ or
$\epsilon\geqslant r-3$.

\begin{lemma}
\label{lemma:surfaces-I} Suppose that $\Sigma\subset\Pi$, where
$\Pi$ is a hyperplane in $\mathbb{P}^{3}$. Then there is a surface
of degree $3r-4-\epsilon$ in $\mathbb{P}^{3}$ that contains the
set $\Sigma\setminus P$ and does not contain the point $P$.
\end{lemma}

\begin{proof}
Suppose that $|\Sigma\setminus
P|>\lfloor(3r-1-\epsilon)/2\rfloor^{2}$. Then
$$
\big(2r-1\big)\big(r-\epsilon\big)-2\geqslant \big|\Sigma\setminus P\big|\geqslant\Big\lfloor\frac{3r-1-\epsilon}{2}\Big\rfloor^{2}+1\geqslant\frac{(3r-2-\epsilon)^{4}}{4}+1,%
$$
which implies that $(r-4)^{2}+2\epsilon r+\epsilon^{2}\leqslant
0$. We have $r=4$ and $\epsilon=0$, which implies that
$$
\big|\Sigma\setminus P\big|\leqslant
\Big\lfloor\frac{3r-1-\epsilon}{2}\Big\rfloor\Big(3r-1-\epsilon-\Big\lfloor\frac{3r-1-\epsilon}{2}\Big\rfloor\Big).
$$

Thus, we proved that in every possible case the inequality
$$
\big|\Sigma\setminus P\big|\leqslant
\mathrm{max}\left(\Big\lfloor\frac{3r-1-\epsilon}{2}\Big\rfloor\Big(3r-1-\epsilon-\Big\lfloor\frac{3r-1-\epsilon}{2}\Big\rfloor\Big),
\Big\lfloor\frac{3r-1-\epsilon}{2}\Big\rfloor^{2}\right).
$$
holds, but at most $3r-4-\epsilon$ points of $\Sigma\setminus P$
can lie on a line, because $3r-4-\epsilon\geqslant 2r-1$.

Let us prove that at most $k(3r-1-\epsilon-k)-2$ points of the set
$\Sigma\setminus P$ can lie on a curve of degree $k\leqslant
(3r-1-\epsilon)/2$. It is enough to show that
$$
k\big(3r-1-\epsilon-k\big)-2\geqslant k\big(2r-1\big)
$$
for all $k\leqslant (3r-1-\epsilon)/2$. We must prove the latter
inequality only for $k>1$ such that
$$
k\big(3r-1-\epsilon-k\big)-2<\big|\Sigma\setminus P\big|\leqslant \big(2r-1\big)\big(r-\epsilon\big)-2,%
$$
because otherwise the condition that at most $k(3r-1-k)-2$ points
of $\Sigma\setminus P$ can lie on a curve of degree $k$ is
vacuous. In particular, we may assume that $k<r-\epsilon$, but
$$
k(3r-1-\epsilon-k)-2\geqslant k(2r-1) \iff r>k-\epsilon,
$$
which implies that at most $k(3r-1-\epsilon-k)-2$ points of
$\Sigma\setminus P$ lie on a curve of degree $k$.

It follows from Theorem~\ref{theorem:Bese} that there is a curve
$C\subset\Pi$ of degree $3r-4-\epsilon$ that contains the set
$\Sigma\setminus P$ and does not contain the point $P$. Let $Y$ be
a sufficiently general cone in $\mathbb{P}^{3}$ over the curve
$C$. Then $Y$ is the required surface.
\end{proof}

Fix a sufficiently general hyperplane $\Pi\subset\mathbb{P}^{3}$.
Let $\psi\colon \mathbb{P}^{3}\dasharrow\Pi$ be a projection from
a sufficiently general point $O\in\mathbb{P}^{3}$. Put
$\Sigma^{\prime}=\psi(\Sigma)$ and $P^{\prime}=\psi(P)$.

\begin{lemma}
\label{lemma:surfaces-II} Suppose that at most $(2r-1)k$ points of
the set $\Sigma^{\prime}$ lie on a possibly reducible curve in
$\Pi$ of degree $k$. Then there is a surface in $\mathbb{P}^{3}$
of degree $3r-4-\epsilon$ that contains all points of the set
$\Sigma\setminus P$ but does not contain the point $P$.
\end{lemma}

\begin{proof}
Arguing as in the proof of Lemma~\ref{lemma:surfaces-I} we obtain
a curve $C\subset\Pi$ of degree $3r-4-\epsilon$ that contains
$\Sigma^{\prime}\setminus P^{\prime}$ and does not pass through
$P^{\prime}$. Let $Y$ be the cone in $\mathbb{P}^{3}$ over the
curve $C$ with the vertex $O$. Then $Y$ is the required surface.
\end{proof}

To conclude the proof of Proposition~\ref{proposition:3r-4} we may
assume that at least $(2r-1)k+1$ points of $\Sigma^{\prime}$ lie
on a curve of degree $k$, where $k$ is the smallest number of such
property.

\begin{lemma}
\label{lemma:conic} The inequality $k\geqslant 3$ holds.
\end{lemma}

\begin{proof}
Let $\Phi\subseteq\Sigma$ be a subset such that $|\Phi|>2(2r-1)$,
but the set $\psi(\Phi)$ is contained in a conic $C\subset\Pi$.
Then the conic $C$ is irreducible. Let $\mathcal{D}$ be a linear
system of quadric surfaces in $\mathbb{P}^{3}$ containing $\Phi$.
Then the base locus of  $\mathcal{D}$ is zero-dimensional by
Lemma~\ref{lemma:zero-dimensional}.

The inequality $k\geqslant 2$ holds by
Lemma~\ref{lemma:zero-dimensional}, which~implies~$r\geqslant 3$.

Let $W$ be a cone in $\mathbb{P}^{3}$ over $C$ with the vertex
$\Omega$. Then
$$
8=D_{1}\cdot D_{2}\cdot W\geqslant
\sum_{\omega\in\Phi}\mathrm{mult}_{\omega}(D_{1})\mathrm{mult}_{\omega}(D_{2})\geqslant |\,\Phi|>2(2r-1)\geqslant 8,%
$$
where $D_{1}$ and $D_{2}$ are general divisors in $\mathcal{D}$,
which is a contradiction.
\end{proof}

There is a subset $\Lambda_{k}^{1}\subseteq\Sigma$ such that
$|\Lambda_{k}^{1}|>(2r-1)k$, but $\psi(\Lambda_{k}^{1})$ is
contained in an irreducible curve of degree $k$. Similarly, we get
a disjoint union
$\cup_{j=k}^{l}\cup_{i=1}^{c_{j}}\Lambda_{j}^{i}$, where
$\Lambda_{j}^{i}$ is a subset in $\Sigma$ such that
$|\Lambda_{j}^{i}|>(2r-1)j$, the points of the subset
$\psi(\Lambda_{j}^{i})$ lie on an irreducible reduced curve in
$\Pi$ of degree $j$, and at most $(2r-1)\zeta$ points of the
subset
$$
\psi\Big(\Sigma\setminus\Big(\bigcup_{j=k}^{l}\bigcup_{i=1}^{c_{j}}\Lambda_{j}^{i}\Big)\Big)\subsetneq\Sigma^{\prime}\subset\Pi\cong\mathbb{P}^{2}
$$
lie on a curve in $\Pi$ of degree $\zeta$. Put
$\Lambda=\cup_{j=k}^{l}\cup_{i=1}^{c_{j}}\Lambda_{j}^{i}$. Let
$\Xi_{j}^{i}$ be the base locus of the linear system of surfaces
of degree $j$ that contains $\Lambda_{j}^{i}$. Then $\Xi_{j}^{i}$
is a finite set by Lemma~\ref{lemma:zero-dimensional} and

\begin{equation}
\label{equation:number-of-good-points}
0\leqslant\big|\Sigma\setminus\Lambda\big|<\big(2r-1\big)\big(r-\epsilon\big)-1-\sum_{i=k}^{l}c_{i}\big(2r-1\big)i<\big(2r-1\big)\Big(r-\epsilon-\sum_{i=k}^{l}ic_{i}\Big).
\end{equation}

\begin{corollary}
\label{corollary:from-the-number-of-good-points} The inequality
$\sum_{i=k}^{l}ic_{i}\leqslant r-\epsilon-1$ holds.
\end{corollary}

We have $\Lambda_{j}^{i}\subseteq\Xi_{j}^{i}$ by construction, but
the points of the set $\Xi_{j}^{i}$ impose independent linear
conditions on homogeneous forms of degree $3(j-1)$ by the
following result.

\begin{lemma}
\label{lemma:non-vanishing} Let $\mathcal{M}$ be a linear
subsystem in $|\mathcal{O}_{\mathbb{P}^n}(\lambda)|$ such that the
base locus of the~li\-near system $\mathcal{M}$ is
zero-dimensional. Then the points of the base locus of
$\mathcal{M}$~impose~inde\-pen\-dent linear conditions on
homogeneous forms of degree $n(\lambda-1)$.
\end{lemma}

\begin{proof}
See Lemma~22 in \cite{Ch04t} or Theorem~3 in \cite{DaGeOr85}.
\end{proof}

Put $\Xi=\cup_{j=k}^{l}\cup_{i=1}^{c_{j}}\Xi_{j}^{i}$. Then
$\Lambda\subseteq\Xi$.

\begin{lemma}
\label{lemma:surfaces-III} Suppose that $\Sigma\subseteq\Xi$. Then
there is a surface in $\mathbb{P}^{3}$ of degree $3r-4-\epsilon$
that contains all points of the set $\Sigma\setminus P$ and does
not contain the point $P\in\Sigma$.
\end{lemma}

\begin{proof}
It follows from Lemma~\ref{lemma:non-vanishing} that for every set
$\Xi_{j}^{i}$ containing the point $P$ there is a surface in
$\mathbb{P}^{3}$ of degree $3(j-1)$ that contains the set
$\Xi_{j}^{i}\setminus P$ and does not contain the point $P$. For
every set $\Xi_{j}^{i}$ not containing the point $P$ there is a
surface of degree $j$ that contains the set $\Xi_{j}^{i}$ and does
not contain $P$ by the definition of the set $\Xi_{j}^{i}$.

The inequality $j<3(j-1)$ holds, because $k\geqslant 2$.
Therefore, for every $\Xi_{j}^{i}\ne\varnothing$ there is a
surface $F_{i}^{j}\subset\mathbb{P}^{3}$ of degree $3(j-1)$ that
contains the set $\Xi_{j}^{i}\setminus(\Xi_{j}^{i}\cap P)$ and
does not contain the point $P$. The union
$\cup_{j=k}^{l}\cup_{i=1}^{c_{j}}F_{j}^{i}$ is a surface of degree
$$
\sum_{i=k}^{l}3(i-1)c_{i}\leqslant\sum_{i=k}^{l}3ic_{i}-3c_{k}\leqslant 3r-6-3\epsilon\leqslant 3r-4-\epsilon%
$$
that contains all points of the set $\Sigma\setminus P$ and does
not contain the point $P$.
\end{proof}

The proof of Lemma~\ref{lemma:surfaces-III} implies that there is
surface of degree $\sum_{i=k}^{l}3(i-1)c_{i}$ that contains
$(\Xi\cap\Sigma)\setminus(\Xi\cap P)$ and does not contain $P$,
and there is a surface of degree $\sum_{i=k}^{l}ic_{i}$ that
contains $\Xi\cap\Sigma$ and does not contain any point of the set
$\Sigma\setminus(\Xi\cap\Sigma)$.

\begin{lemma}
\label{lemma:swapping} Let $\Lambda$ and $\Delta$ be disjoint
finite subsets in $\mathbb{P}^{n}$ such that there is a
hypersurface of degree $\zeta\leqslant\xi$ that contains $\Lambda$
and does not contain any point in $\Delta$, the points of
$\Lambda$~impose independent linear conditions on hypersurfaces of
degree $\xi$, the points of~$\Delta$~impose~in\-de\-pen\-dent
linear conditions on hypersurfaces of degree $\xi-\zeta$. Then the
points of~$\Lambda\cup\Delta$~impose independent linear conditions
on hypersurfaces of degree $\xi$.
\end{lemma}

\begin{proof}
Let $Q$ be a point in $\Lambda\cup\Delta$. To conclude the proof
we must find a hypersurface of degree $\xi$ that contains
$(\Lambda\cup\Delta)\setminus Q$ and does not contain $Q$. We may
assume that $Q\in\Lambda$.

Let $F$ be the homogenous form of degree $\xi$ that vanishes at
$\Lambda\setminus Q$ and does not vanish at $Q$. Put
$\Delta=\{Q_{1},\ldots,Q_{\delta}\}$, where $Q_{i}$ is a point.
Then there~is~a~homo\-ge\-neous form $G_{i}$ of degree $\xi$ that
vanishes at $(\Lambda\cup\Delta)\setminus Q_{i}$ and does not
vanish at $Q_{i}$. We have
$$
F\big(Q_{i}\big)+\mu_{i}G_{i}\big(Q_{i}\big)=0
$$
for some $\mu_{i}\in\mathbb{C}$, because $g_{i}(Q_{i})\ne 0$. Then
the homogenous form $F+\sum_{i=1}^{\delta}\mu_{i}G_{i}$ vanishes
at every point of the set $(\Lambda\cup\Delta)\setminus Q$ and
does not vanish at the point $Q$.
\end{proof}

Put $d=3r-4-\epsilon-\sum_{i=k}^{l}ic_{i}$ and
$\bar{\Sigma}=\psi(\Sigma\setminus(\Xi\cap\Sigma))$. It follows
from Lemma~\ref{lemma:swapping} that to prove
Proposition~\ref{proposition:3r-4-epsilon} it is enough to show
that the points of the subset $\bar{\Sigma}\subset\Pi$ and the
integer~$d$ satisfy the hypotheses of Theorem~\ref{theorem:Bese}.
We may assume that
$\emptyset\ne\bar{\Sigma}\subsetneq\Sigma^{\prime}$.

\begin{lemma}
\label{lemma:surfaces-IV} The inequality $|\bar{\Sigma}|\leqslant
\lfloor (d+3)/2\rfloor^{2}$ holds.
\end{lemma}

\begin{proof}
Suppose that the inequality $|\bar{\Sigma}|\geqslant\lfloor
(d+3)/2\rfloor^{2}+1$ holds. Then
$$
\big(2r-1\big)\big(r-\epsilon-\sum_{i=k}^{l}c_{i}i\big)-2\geqslant\big|\bar{\Sigma}\big|\geqslant \Big\lfloor\frac{d+2}{4}\Big\rfloor^{2}+1\geqslant\frac{(3r-2-\epsilon-\sum_{i=k}^{l}ic_{i})^{2}}{4}+1%
$$
by Corollary~\ref{corollary:from-the-number-of-good-points}. Put
$\Delta=\epsilon+\sum_{i=k}^{l}c_{i}i$. Then $\Delta\geqslant
k\geqslant 3$ and
$$
4\big(2r-1\big)\big(r-\Delta\big)-12\geqslant\big(3r-2-\Delta\big)^{2},
$$
which implies that $r^{2}-8r+16+2r\Delta+\Delta^{2}\leqslant 0$,
which is a contradiction.
\end{proof}

The inequality $d\geqslant 3$ holds by
Corollary~\ref{corollary:from-the-number-of-good-points}, because
$r\geqslant 3$.

\begin{lemma}
\label{lemma:surfaces-V} Suppose that at least $d+1$ points of the
set $\bar{\Sigma}$ lie on a line. Then there is a surface in
$\mathbb{P}^{3}$ of degree $3r-4-\epsilon$ containing
$\Sigma\setminus P$ and not passing through the point $P$.
\end{lemma}

\begin{proof}
We have $|\bar{\Sigma}|\geqslant d+1$. Hence, it follows from the
inequalities~\ref{equation:number-of-good-points} that
$$
3r-3-\epsilon-\sum_{i=k}^{l}ic_{i}<\big(2r-1\big)\big(r-\epsilon\big)-1-\sum_{i=k}^{l}c_{i}\big(2r-1\big)i,%
$$
which gives $\sum_{i=k}^{l}ic_{i}\ne r-\epsilon-1$. Now it follows
from Corollary~\ref{corollary:from-the-number-of-good-points} that
$\sum_{i=k}^{l}ic_{i}\leqslant r-\epsilon-2$, but $2r-1\geqslant
3r-3-\epsilon-\sum_{i=k}^{l}ic_{i}$, which implies that
$\sum_{i=k}^{l}ic_{i}=r-\epsilon-2$ and $d=2r-2$.

We have a surface of degree $\sum_{i=k}^{l}3(i-1)c_{i}\leqslant
3r-4-\epsilon$ that contains $(\Xi\cap\Sigma)\setminus (\Xi\cap
P)$~and does not contain the point $P$, and we have a surface of
degree $r-\epsilon-2$ that contains all points of the set
$\Xi\cap\Sigma$ and does not contain any point of the set
$\Sigma\setminus(\Xi\cap\Sigma)$.

The set $\Sigma\setminus(\Xi\cap\Sigma)$ contains at most $4r-4$
points, but at most $2r-1$ points of $\Sigma$ lie on a line. The
points of $\Sigma\setminus(\Xi\cap\Sigma)$ impose independent
linear conditions on homo\-ge\-ne\-ous forms of degree $2r-2$ by
Theorem~2 in \cite{EiJ87}, which concludes the proof by
Lemma~\ref{lemma:swapping}.
\end{proof}

Therefore, we may assume that at most $d$ points of the set
$\bar{\Sigma}$ lie on a line in $\Pi$.

\begin{lemma}
\label{lemma:surfaces-VI} At most $t(d+3-t)-2$ points of
$\bar{\Sigma}$ lie on a curve in $\Pi$ of degree $t\leqslant
(d+3)/2$.
\end{lemma}

\begin{proof}
At most $(2r-1)t$ of the points of $\bar{\Sigma}$ lie on a curve
in $\Pi$ of degree $t$, which implies that to conclude the proof
it is enough to show that the inequality
$$
t\big(d+3-t\big)-2\geqslant \big(2r-1\big)t
$$
holds for every $t\leqslant (d+3)/2$ such that $t>1$ and
$t(d+3-t)-2<|\bar{\Sigma}|$. We have
$$
t(d+3-t)-2\geqslant t(2r-1)\iff
t(r-\epsilon-\sum_{i=k}^{l}ic_{i}-t)\geqslant 2\iff
r-\epsilon-\sum_{i=k}^{l}ic_{i}>t,
$$
because $t>1$. Therefore, we may assume that the inequalities
$t(d+3-t)-2<|\bar{\Sigma}|$ and
$$
r-\epsilon-\sum_{i=k}^{l}ic_{i}\leqslant t\leqslant {\frac{d+3}{2}}%
$$
hold. Let $g(x)=x(d+3-x)-2$. Then $g(x)$ is increasing for every
$x<(d+3)/2$, which implies that $g(t)\geqslant
g(r-\epsilon-\sum_{i=k}^{l}ic_{i})$. Now the
inequalities~\ref{equation:number-of-good-points} imply that
$$
\big(2r-1\big)\Big(r-\epsilon-\sum_{i=k}^{l}ic_{i}\Big)-2\geqslant\big|\bar{\Sigma}\big|>g(t)\geqslant \Big(r-\epsilon-\sum_{i=k}^{l}ic_{i}\Big)\big(2r-1\big)-2,%
$$
which is a contradiction.
\end{proof}

We can apply Theorem~\ref{theorem:Bese} to the blow up of $\Pi$ at
the points of $\bar{\Sigma}$ and the integer $d$, which implies
the existence of a surface in $\mathbb{P}^{3}$ of degree
$3r-4-\epsilon$ that contains every point of the set
$\Sigma\setminus P$ and does not contain the point $P$ by
Lemma~\ref{lemma:swapping}.

\section{Auxiliary result.}
\label{section:extra}

In this section we prove Theorem~\ref{theorem:Hong-and-Park}. Let
$\pi\colon X\to\mathbb{P}^{3}$ be a double cover branched over a
surface $S$ of degree $2r\geqslant 4$ with isolated ordinary
double points.

\begin{lemma}
\label{lemma:star-property} Let $F$ be a hypersurface in
$\mathbb{P}^{n}$ of degree $d$ such that $F$ has isolated
sin\-gu\-la\-ri\-ties, and $C$ be a curve in $\mathbb{P}^{n}$ of
degree $k$. Then $C$ contains at most $k(d-1)$ singular points of
the hypersurface $F$, and the equality
$|\,\mathrm{Supp}(C)\cap\mathrm{Sing}(F)|=k(d-1)$ implies that
every singular point of the hypersurface $F$ contained in $C$ is
non-singular on  the curve $C$.
\end{lemma}

\begin{proof}
Let $f(x_0,\ldots,x_n)$ be the homogeneous form of degree $d$ such
that $f(x_0,\ldots,x_n)=0$ defines the hypersurface $F$, where
$(x_0:\ldots:x_n)$ are homogeneous coordinates on
$\mathbb{P}^{n}$. Put
$$
\mathcal{D}=\left|\sum_{i=0}^{n}\lambda_i\frac{\partial f}{\partial x_i}=0\right|\subset\big|\mathcal{O}_{\mathbb{P}^n}(d-1)\big|,%
$$
where $\lambda_{i}\in\mathbb{C}$. Then the base locus of the
linear system $\mathcal{D}$ consists of singular points of the
hypersurface $F$. Therefore, the curve $C$ intersects a generic
member of the linear system $\mathcal{D}$ at most $(d-1)k$ times,
which implies the claim.
\end{proof}

\begin{lemma}
\label{lemma:Hong-Park-01} Suppose that there are plane
$\Pi\subset\mathbb{P}^{3}$ and a reduced curve $C\subset\Pi$ of
degree~$r$ that contains $(2r-1)r$ singular points of $S$. Then
$S$ can be defined by the equation~\ref{equation:nodal-surface}.
\end{lemma}

\begin{proof}
Let $S\vert_{\Pi}=\sum_{i=1}^{\alpha}m_{i}C_{i}$, where $C_{i}$ is
an irreducible reduced curve, and $m_{i}$ is a natural number. We
may assume that $C_{i}\ne C_{j}$ for $i\ne j$, and
$C=\sum_{i=1}^{\beta}C_{i}$, where $\beta\leqslant\alpha$. Then
\begin{equation}
\label{equation:degrees}
\sum_{i=1}^{\beta}\mathrm{deg}(C_{i})=r=\frac{\sum_{i=1}^{\alpha}m_{i}\mathrm{deg}(C_{i})}{2},
\end{equation}
which implies that the curve $C_{i}$ contains exactly
$(2r-1)\mathrm{deg}(C_{i})$ singular points of the surface $S$ for
every $i\leqslant\beta$ due to Lemma~\ref{lemma:star-property}.
Moreover, the curve $C$ is smooth at every singular point of the
surface $S$ that is contained in the curve $C$ by
Lemma~\ref{lemma:star-property}.

Suppose that $m_{\gamma}=1$ for some $\gamma\leqslant\beta$. Then
$C_{\gamma}$ contains $(2r-1)\mathrm{deg}(C_{\gamma})$ singular
points of the surface $S$, but the curve $S\vert_{\Pi}$ must be
singular at every singular point of the surface $S$ that is
contained in $C_{\gamma}$. Thus, we have
$$
\mathrm{Sing}\big(S\big)\cap\mathrm{Supp}\big(C_{\gamma}\big)\subseteq\bigcup_{i\ne\gamma}C_{i}\cap C_{\gamma},%
$$
but $|\,C_{i}\cap C_{\gamma}|\leqslant (C_{i}\cdot
C_{\gamma})_{\Pi}=\mathrm{deg}(C_{i})\mathrm{deg}(C_{\gamma})$ for
$i\ne\gamma$. Hence, we have
$$
\sum_{i\ne\gamma}\mathrm{deg}(C_{i})\mathrm{deg}(C_{\gamma})\geqslant \big(2r-1\big)\mathrm{deg}(C_{\gamma}),%
$$
but on the plane $\Pi$ we have the equalities
$$
\big(2r-\mathrm{deg}(C_{\gamma})\big)\mathrm{deg}(C_{\gamma})=\big(S\vert_{\Pi}-C_{\gamma}\big)\cdot C_{\gamma}=\sum_{i\ne\gamma}m_{i}\mathrm{deg}(C_{i})\mathrm{deg}(C_{\gamma}),%
$$
which implies that $\mathrm{deg}(C_{\gamma})=1$ and $m_{i}=1$ for
every $i$. Now the equalities~\ref{equation:degrees} imply that
the equality $\beta<\alpha$ holds, but every singular point of the
surface $S$ that is contained in the curve $C$ must be an
intersection point of $C$ and the curve
$\sum_{i=\beta+1}^{\alpha}C_{i}$, which consists of at most
$r^{2}$ points, which is a contradiction.

Hence, we have $m_{i}\geqslant 2$ for every $i\leqslant\beta$.
Therefore, it follows from the equalities~\ref{equation:degrees}
that $\alpha=\beta$ and $m_{i}=2$ for every $i$.

Let $f(x,y,z,w)$ be the homogeneous form of degree $2r$ such that
$f=0$ defines the surface $S$, where $(x:y:z:w)$ are homogeneous
coordinates on $\mathbb{P}^{3}$. We may assume that the plane
$\Pi$ is given by the equation $x=0$. Then
$f(0,y,z,w)=g_{r}^{2}(y,z,w)$, where $g_{r}$ is a homogeneous
polynomial of degree $r$ such that $C$ is given by $x=g_{r}=0$,
which implies that the surface $S$ can be defined by the
equation~\ref{equation:nodal-surface}.
\end{proof}

It follows from Lemma~\ref{lemma:star-property} that at most
$(2r-1)k$ singular points of the surface $S$ can lie on a curve of
degree $k$. However, the claim of Lemma~\ref{lemma:star-property}
can be improved for curves that are not contained in
two-dimensional linear subspaces of $\mathbb{P}^{3}$.

\begin{lemma}
\label{lemma:Hong-Park-02} Let $C$ be an irreducible reduced curve
in $\mathbb{P}^{3}$ of degree $k$ that is not contained in a
hyperplane. Then $|C\cap\mathrm{Sing}(S)|\leqslant(2r-1)k-2$.
\end{lemma}

\begin{proof}
Suppose that the curve $C$ contains at least $(2r-1)k-1$ singular
points of the surface $S$. Then $C\subset S$, because otherwise we
have
$$
2rk=\mathrm{deg}(C)\mathrm{deg}(S)\leqslant 2(2r-1)k-2=4rk-2k-2,
$$
which leads to $2k(r-1)\leqslant 2$, but $r\geqslant 2$ and
$k\geqslant 3$.

Let $O$ be a sufficiently general point of the curve $C$, and
$\psi\colon \mathbb{P}^{3}\dasharrow\Pi$ be a projection from the
point $O$, where $\Pi$ is a sufficiently general plane in
$\mathbb{P}^{3}$. Then $\psi\vert_{C}$ is a birational mor\-phism,
because $C$ is not a plane curve. Put $Z=\psi(C)$. Then $Z$ has
degree $k-1$.

Let $Y$ be a cone in $\mathbb{P}^{3}$ over  $Z$ with the vertex
$O$. Then $C\subset Y$.

It follows from the generality of the point $O$ that the point $O$
is not contained in a hyperplane in $\mathbb{P}^{3}$ that is
tangent to the surface $S$ at some point of the curve $C$, because
the curve $C$ is not contained in a hyperplane. Therefore, the
cone $Y$ does not tangent the surface $S$ along the curve $C$.

Put $S\vert_{Y}=C+R$, where $R$ is a curve of degree $2rk-k-2r$.
Then the generality of the point $O$ implies that the curve $R$
does not contains rulings of the cone $Y$.

Let $\alpha:\bar{Z}\to Z$ be the normalization of $Z$. Then there
is a commutative diagram
$$
\xymatrix{
\bar{Y}\ar@{->}[rr]^{\beta}\ar@{->}[d]_{\pi}&&Y\ar@{-->}[d]^{\psi\vert_{Y}}\\%
\bar{Z}\ar@{->}[rr]_{\alpha}&&Z,}
$$ %
where $\beta$ is a birational morphism, $\bar{Y}$ is smooth, and
$\pi$ is a $\mathbb{P}^{1}$-bundle. Let $L$ be a general fiber of
$\pi$, and $E$ be a section of $\pi$ such that $\beta(E)=O$. Then
$E^{2}=-k+1$ on $\bar{Y}$.

Let $\bar{C}$ and $\bar{R}$ be proper transforms of the curves $C$
and $R$ on the surface $\bar{Y}$ respectively, and $Q$ be an
arbitrary point of the set $\mathrm{Sing}(S)\cap C$. Then there is
a point $\bar{Q}\in\bar{Y}$ such that $\beta(\bar{Q})=Q$ and
$\bar{Q}\in\mathrm{Supp}(\bar{C}\cdot \bar{R})$, but
$$
\bar{R}\equiv (2r-2)E+(2rk-k-2r)L
$$
and $\bar{C}\equiv E+kL$. Therefore, we have
$(2r-1)k-2=\bar{C}\cdot\bar{R}\geqslant (2r-1)k-1$.
\end{proof}

Now we prove Theorem~\ref{theorem:Hong-and-Park} by reductio ad
absurdum. Put $\Sigma=\mathrm{Sing}(S)$, and suppose that the
following~conditions~hold:
\begin{itemize}
\item the inequalities $|\Sigma|\leqslant(2r-1)r+1$ and $r\geqslant 3$ hold;%
\item the surface $S$ can not be defined by the equation~\ref{equation:nodal-surface};%
\item the threefold $X$ is not factorial, which implies that there is a point $P\in\Sigma$~such~that every surface in $\mathbb{P}^{3}$ of degree $3r-4$ containing $\Sigma\setminus P$ contains the point $P$.%
\end{itemize}

We assume that $r\geqslant 4$, because the case $r=3$ is done in
\cite{HoPa04}.

\begin{lemma}
\label{lemma:Hong-Park-03} Let $\Pi$ be a two-dimensional linear
subspace in $\mathbb{P}^{3}$. Then
$\vert\Pi\cap\Sigma\vert\leqslant 2r$.
\end{lemma}

\begin{proof}
Suppose that $\vert\Pi\cap\Sigma\vert>2r$. Let us show that this
assumption leads to a contradiction. Let $\Gamma$ be the subset of
the set $\Sigma$ that consists of all points that are not
contained in the plane $\Pi$. Then $\Gamma$ contains at most
$(2r-1)(r-1)-1$ points, which impose independent linear conditions
on homogeneous forms of degree $3r-5$ by
Proposition~\ref{proposition:3r-4-epsilon}.

Suppose that $P\not\in\Pi$. Then there is a surface
$F\subset\mathbb{P}^{3}$ of degree $3r-5$ that contains the set
$\Gamma\setminus P$ and does not contain the point $P$. Hence, the
union $F\cup\Pi$ is the surface of degree $3r-4$ that contains the
set $\Sigma\setminus P$ and does not contain the point $P$, which
is impossible due to our assumptions. Therefore, we have
$P\in\Pi$.

The curve $\Pi\cap S$ is singular in every point of the set
$\Pi\cap\Sigma$. Thus, it follows from the proof of
Lemma~\ref{lemma:star-property} that
$|\Pi\cap\Sigma|\leqslant(2r-1)r$, but
Lemma~\ref{lemma:Hong-Park-01} implies that $\Pi\cap\Sigma$ is not
contained in a curve of degree $r$ if $|\Pi\cap\Sigma|=(2r-1)r$.
The proof of Lemma~\ref{lemma:surfaces-I} implies that there is a
surface of degree $3r-4$ that contains the set
$(\Pi\cap\Sigma)\setminus P$ and does not contain the point $P$,
which concludes the proof by Lemma~\ref{lemma:swapping}.
\end{proof}

The inequality $|\Sigma|\geqslant(2r-1)r$ holds by
Proposition~\ref{proposition:3r-4}.

\begin{lemma}
\label{lemma:Hong-Park-04} Let $L_{1}$ and $L_{2}$ be distinct
lines in $\mathbb{P}^{3}$. Then $|(L_{1}\cup
L_{2})\cap\Sigma|<4r-2$.
\end{lemma}

\begin{proof}
Suppose that $|(L_{1}\cup L_{2})\cap\Sigma|\geqslant 4r-2$. Then
$|L_{i}\cap\Sigma|=2r-1$ by Lemma~\ref{lemma:star-property}, and
the lines $L_{1}$ and $L_{2}$ are not contained in one hyperplane
by Lemma~\ref{lemma:Hong-Park-03}.

Fix two points $Q_{1}$ and $Q_{2}$ in $\Sigma\setminus((L_{1}\cup
L_{2})\cap\Sigma)$ different from $P$ such that $Q_{1}\ne Q_{2}$,
and let $\Pi_{i}$ be a plane in $\mathbb{P}^{3}$ that contains
$L_{i}$ and $Q_{i}$. Then $|\Pi_{i}\cap\Sigma|=2r$ by
Lemma~\ref{lemma:Hong-Park-03}.

Suppose that $P\not\in\Pi_{1}\cup\Pi_{2}$. Then there is a surface
$F\subset\mathbb{P}^{3}$ of degree $3r-6$ that does not contain
the point $P$ and contains all points of the set
$$
\Big(\Sigma\setminus\Big(\Sigma\cap\big(\Pi_{1}\cup\Pi_{2}\big)\Big)\Big)\setminus P%
$$
by Proposition~\ref{proposition:3r-4-epsilon}. Hence, the surface
$F\cup\Pi_{1}\cup\Pi_{2}$ is a surface in $\mathbb{P}^{3}$ of
degree $3r-4$ that contains all points of the set $\Sigma\setminus
P$ and does not contain the point $P$, which contradicts to our
assumption. Therefore, we have $P\in\Pi_{1}\cup\Pi_{2}$.

The set $\Sigma\cap(\Pi_{1}\cup\Pi_{2})$ consists of $4r$ points
by Lemma~\ref{lemma:Hong-Park-03}. Therefore, the points of the
set $\Sigma\cap(\Pi_{1}\cup\Pi_{2})$ impose independent linear
conditions on homogeneous forms $\mathbb{P}^{3}$ of degree $3r-4$
by Theorem~2 in \cite{EiJ87}. On the other hand, the inequality
$$
\Big|\Sigma\setminus\Big(\Sigma\cap\big(\Pi_{1}\cup\Pi_{2}\big)\Big)\Big|<\big(2r-1\big)\big(r-2\big)
$$
holds, and the points of
$\Sigma\setminus(\Sigma\cap(\Pi_{1}\cup\Pi_{2}))$ impose
independent linear conditions homogeneous forms of degree $3r-6$
by Proposition~\ref{proposition:3r-4-epsilon}, which is impossible
by Lemma~\ref{lemma:swapping}.
\end{proof}

\begin{lemma}
\label{lemma:Hong-Park-05} Let $C$ be a curve in $\mathbb{P}^{3}$
of degree $k\geqslant 2$. Then $|C\cap\Sigma|<(2r-1)k$.
\end{lemma}

\begin{proof}
Suppose that $|C\cap\Sigma|\geqslant (2r-1)k$. Let us show that
this assumption leads to a contradiction. We have
$|C\cap\Sigma|=(2r-1)k$ by Lemma~\ref{lemma:star-property}, and
the curve $C$ is not contained in a hyperplane by
Lemma~\ref{lemma:Hong-Park-03}. Therefore, the curve $C$ is
reducible by Lemma~\ref{lemma:Hong-Park-02}.

Let us put $C=\sum_{i=1}^{\alpha}C_{i}$, where $\alpha\geqslant 2$
and $C_{i}$ is an irreducible curve. Then
$k=\sum_{i=1}^{\alpha}d_{i}$, where $d_{i}$ is the degree of the
curve $C_{i}$, which implies $|C_{i}\cap\Sigma|=(2r-1)d_{i}$ by
Lemma~\ref{lemma:star-property}.

The curve $C_{i}$ is contained in a hyperplane in $\mathbb{P}^{3}$
by Lemma~\ref{lemma:Hong-Park-02}. So, the equalities $d_{i}=1$
and $\alpha=k$ hold by Lemma~\ref{lemma:Hong-Park-03} for all $i$,
which contradicts Lemma~\ref{lemma:Hong-Park-04}, because
$k\geqslant 2$.
\end{proof}

\begin{lemma}
\label{lemma:Hong-Park-06} Let $L$ be a line in $\mathbb{P}^{3}$.
Then $|L\cap\Sigma|\leqslant 2r-2$.
\end{lemma}

\begin{proof}
Suppose that the inequality $|L\cap\Sigma|\geqslant 2r-1$ holds.
Let us show that this assumption leads to a contradiction. We have
$|L\cap\Sigma|=2r-1$ by Lemma~\ref{lemma:star-property}.

Let $\Phi$ be a hyperplane in $\mathbb{P}^{3}$ such that $\Phi$
contains the line $L$, and $\Phi$ contains an arbitrary point of
the set $\Sigma\setminus(L\cap\Sigma)$. Then $\Phi$ contains $2r$
points of the set $\Sigma$ by Lemma~\ref{lemma:Hong-Park-03}.

Put $\Delta=\Sigma\setminus(\Phi\cap\Sigma)$. Then
$|\Delta|\leqslant (2r-1)(r-1)$.

Suppose that the points of the set $\Delta$ impose independent
linear conditions on homogeneous forms on $\mathbb{P}^{3}$ of
degree $3r-5$. Then it follows from Lemma~\ref{lemma:swapping}
that the points of the set $\Sigma$ impose independent linear
conditions on homogeneous forms of degree $3r-4$, because the
points of the set $\Phi\cap\Sigma$ impose independent linear
conditions on homogeneous forms on $\mathbb{P}^{3}$ of degree
$3r-4$. Therefore, the points of the set $\Delta$ impose dependent
linear conditions on homogeneous forms on $\mathbb{P}^{3}$ of
degree $3r-5$.

There is a point $Q\in\Delta$ such that every surface of degree
$3r-5$ containing $\Delta\setminus Q$~must contain $Q$, which
implies $|\Delta|=(2r-1)(r-1)$ and $|\Sigma|=(2r-1)r+1$ by
Proposition~\ref{proposition:3r-4-epsilon}.

Fix sufficiently general hyperplane $\Pi\subset\mathbb{P}^{3}$ and
a point $O\in\mathbb{P}^{3}$. Let $\psi\colon
\mathbb{P}^{3}\dasharrow\Pi$ be a projection from the point $O$.
Put $\Delta^{\prime}=\psi(\Delta)$ and $Q^{\prime}=\psi(Q)$. Then
at most $2r-2$ points of the set $\Delta^{\prime}$ lie on a line
by Lemmas~\ref{lemma:zero-dimensional} and
\ref{lemma:Hong-Park-04}.

Suppose that at most $(2r-1)k$ points of the set $\Delta^{\prime}$
lie on any curve in $\Pi$ of degree~$k$ for every natural number
$k$, and there is a curve $Z\subset\Pi$ of degree $r-1$ that
contains the whole set $\Delta^{\prime}$. Then the points of the
set $\Delta$ impose independent linear conditions on homogeneous
forms on $\mathbb{P}^{3}$ of degree $3r-5$ by
Lemmas~\ref{lemma:zero-dimensional},~\ref{lemma:non-vanishing}~and~\ref{lemma:Hong-Park-05}
in the case when the curve $Z$ is irreducible. So, we have
$Z=\sum_{i=1}^{\alpha}Z_{i}$, where $\alpha\geqslant 2$, and
$Z_{i}$~is~an~irre\-du\-cible curve of degree $d_{i}$. Then
$r=\sum_{i=1}^{\alpha}d_{i}$, which implies that $Z_{i}$ contains
$(2r-1)d_{i}$ points of the set $\Delta^{\prime}$, and every point
of the set $\Delta^{\prime}$ is contained in one irreducible
component of the curve $Z$. In particular, we have $d_{i}\ne 1$
for every $i$.

Let $Z_{\beta}$ be the component of $Z$ containing $Q^{\prime}$,
and $\Gamma$ be a subset of $\Delta$ such that
$$
\psi\big(\Gamma\big)=\Delta^{\prime}\cap Z_{\beta}\subset\Pi\cong\mathbb{P}^{2},%
$$
which implies $Q\in\Gamma$. There is a surface
$F_{\beta}\subset\mathbb{P}^{3}$ of degree $3(d_{\beta}-1)$ that
contains all point of the set $\Gamma\setminus Q$~and does not
contain $Q$ by Lemmas~\ref{lemma:zero-dimensional},
\ref{lemma:non-vanishing} and \ref{lemma:Hong-Park-05}. Let
$Y_{i}$~be~a~cone over the curve $Z_{i}$, whose vertex is the
point $O$. Then the union $F_{\beta}\cup\cup_{i\ne\beta}Y_{i}$ is
a surface of degree $3d_{i}-3+\sum_{i\ne\beta}d_{i}=2d_{i}+r-4$
that contains $\Delta\setminus Q$ and does not contains $Q$, which
is impossible, because $2d_{i}+r-4\leqslant 3r-5$. Hence, we
proved that
\begin{itemize}
\item either at least $(2r-1)k+1$ points of $\Delta^{\prime}$ lie on a curve in $\Pi$ of degree $k$;%
\item or there is no curve in $\Pi$ of degree $r-1$ that contains the whole set $\Delta^{\prime}$.%
\end{itemize}

Suppose that at most $(2r-1)k$ points of the set $\Delta^{\prime}$
lie on every curve in $\Pi$ of degree~$k$ for every natural $k$.
Then the points of the set $\Delta^{\prime}\setminus Q^{\prime}$
and the number $3r-5$ satisfy all hypotheses of
Theorem~\ref{theorem:Bese}, because there is no curve in $\Pi$ of
degree $r-1$ that contains the set $\Delta^{\prime}$. Hence, we
can apply Theorem~\ref{theorem:Bese} to the blow up of the plane
$\Pi$ at the points of the set $\Delta^{\prime}\setminus
Q^{\prime}$ to prove the existence of a curve in the plane $\Pi$
of degree $3r-5$ that contains the set $\Delta^{\prime}\setminus
Q^{\prime}$ and does not contains the point $Q^{\prime}$, which is
a contradiction.

Therefore, at least $(2r-1)k+1$ points of the set
$\Delta^{\prime}$ lie on a curve in $\Pi$ of degree $k$, where
$k\geqslant 3$ by Lemma~\ref{lemma:conic}. Thus, the proof of
Proposition~\ref{proposition:3r-4-epsilon} implies the existence
of a subset $\Xi\subseteq\Delta$ such that the following
conditions hold:
\begin{itemize}
\item the points of $\Xi$ impose independent linear conditions on surfaces of degree $3r-5$;%
\item at most $(2r-1)k$ points of the set $\psi(\Delta\setminus\Xi)$ lie on a curve in $\Pi$ of degree $k$;%
\item there is a surface of degree $\mu\leqslant r-2$ that contains all points of the set $\Xi$ and does not contain any point of the set $\Delta\setminus\Xi$;%
\item the inequality $|\Delta\setminus\Xi|\leqslant(2r-1)(r-1-\mu)-1$ holds.%
\end{itemize}

Put $\bar{\Delta}=\psi(\Delta\setminus\Xi)$ and $d=3r-5-\mu$. Then
the points of $\bar{\Delta}$ impose dependent linear conditions on
homogeneous forms of degree $d$ by Lemma~\ref{lemma:swapping},
which implies that there is a point $\bar{Q}\in\bar{\Delta}$ such
that $\bar{\Delta}\setminus\bar{Q}$ and $d$ do not satisfy the
hypotheses of Theorem~\ref{theorem:Bese}.

We have $d\geqslant 3$, because $r\geqslant 4$. The proof of
Lemma~\ref{lemma:surfaces-IV}~gives
$$
\big|\bar{\Delta}\setminus\bar{Q}\big|\leqslant\Big\lfloor\frac{d+3}{2}\Big\rfloor^{2},%
$$
which implies that at least $t(d+3-t)-1$ points of the finite set
$\bar{\Delta}\setminus\bar{Q}$ lie on a curve of degree $t$ for
some natural number $t$ such that $t\leqslant (d+3)/2$.

Suppose that $t=1$. At least $d+1$ points of $\bar{\Delta}$ lie on
a line, but at most $2r-2$ points of $\Delta^{\prime}$ lie on a
line by Lemmas~\ref{lemma:zero-dimensional} and
\ref{lemma:Hong-Park-04}, which implies that $d=2r-3$ and
$|\bar{\Delta}|=2r-2$ points, which is impossible because the
points of the set $\bar{\Delta}$ impose dependent linear
conditions on homogeneous forms of degree~$d$. Therefore, we see
that $t\geqslant 2$.

At least $t(d+3-t)-1$ points of $\bar{\Delta}\setminus\bar{Q}$ lie
on a curve of degree $t\geqslant 2$. Then
$$
t\big(d+3-t\big)-1\leqslant\big|\bar{\Delta}\setminus\bar{Q}\big|
\leqslant\big(2r-1\big)\big(r-1\big)-2-\mu\big(2r-1\big)i,
$$
but $t(d+3-t)-1\leqslant(2r-1)t$, because at most $(2r-1)t$ points
of $\bar{\Delta}$ lie on a curve of degree $t$. Hence, we have
$t\geqslant r-1-\mu$, which gives
$$
\big(2r-1\big)\Big(r-1-\mu\Big)-2\geqslant\big|\bar{\Delta}\setminus\bar{Q}\big|\geqslant t(d+3-t)-1\geqslant \Big(r-1-\mu\Big)\big(2r-1\big)-1,%
$$
which is a contradiction.
\end{proof}

\begin{corollary}
\label{corollary:Hong-Park-strick-inequality-for-curves} Let $C$
be any curve in $\mathbb{P}^{3}$ of degree $k$. Then
$|C\cap\Sigma|<(2r-1)k$.
\end{corollary}

Fix a hyperplane $\Pi\subset\mathbb{P}^{3}$ and a general point
$O\in\mathbb{P}^{3}$. Let
$$
\psi\colon\mathbb{P}^{3}\dasharrow\Pi\subset\mathbb{P}^{3}
$$
be a projection from $O$. Put $\Sigma^{\prime}=\psi(\Sigma)$ and
$P^{\prime}=\psi(P)$. Then
$\psi\vert_{\Sigma}\colon\Sigma\to\Sigma^{\prime}$ is a bijection.

\begin{lemma}
\label{lemma:Hong-Park-07} Let $C$ be an irreducible curve in
$\Pi$ of degree $r$. Then $|\,C\cap\Sigma^{\prime}|<(2r-1)r$.
\end{lemma}

\begin{proof}
Suppose that $|C\cap\Sigma^{\prime}|\geqslant (2r-1)r$. Let us
show that this assumption leads to a contradiction. Let $\Psi$ be
a subset of the set $\Sigma$ consisting of the points that are
mapped to the curve $C$ by the projection $\psi$. Then
$|\Psi|\geqslant(2r-1)r$, but less than $(2r-1)r$ points of the
set $\Sigma$ lie on a curve of degree $r$ by
Corollary~\ref{corollary:Hong-Park-strick-inequality-for-curves}.

Let $\mathcal{H}$ be a linear system of surfaces in
$\mathbb{P}^{3}$ of degree $r$ that contains $\Psi$, and $\Phi$ be
the base locus of $\mathcal{H}$. Then $\Phi$ is finite
Lemma~\ref{lemma:zero-dimensional}. Put $\Upsilon=\Sigma\cap\Phi$.
Then the points of $\Upsilon$ impose independent linear conditions
on homogeneous forms of degree $3r-3$ by
Lemma~\ref{lemma:non-vanishing}.

Let $\Gamma$ be a subset in $\Upsilon$ such that
$\Upsilon\setminus\Gamma$ consists of $4r-6$ points. Then
$$
\big|\Gamma\big|\leqslant 2r^{2}-5r-5\leqslant \frac{(r+2)(r+1)r}{6}-1,%
$$
because $r\geqslant 4$, which implies that there is a surface
$F\subset\mathbb{P}^{3}$ of degree $r-1$ that contains all points
of the set $\Gamma$. Let $\Theta$ be a subset of the set
$\Upsilon$ such that $\Theta$ consists of all points that are
contained in the surface $F$. Then the points of the set $\Theta$
impose independent linear conditions on homogeneous forms on
$\mathbb{P}^{3}$ of degree $3r-4$ by Theorem~3 in \cite{DaGeOr85}.

Put $\Delta=\Upsilon\setminus\Theta$. Then the points of $\Delta$
impose independent linear conditions on homogeneous forms of
degree $2r-3$ by Theorem~2 in \cite{EiJ87} and
Lemmas~\ref{lemma:Hong-Park-03} and \ref{lemma:Hong-Park-06}. So,
the points of the set $\Upsilon$ impose independent linear
conditions on homogeneous forms of degree $3r-4$ by
Lemma~\ref{lemma:swapping}, which also follows from Theorem~3 in
\cite{DaGeOr85}, because $(2r-1)r+1<r^{3}$.

We have $|\Sigma\setminus\Upsilon|\leqslant 1$. Thus, the points
of the set $\Sigma$ impose independent linear conditions on
homogeneous forms of degree $3r-4$ by Lemma~\ref{lemma:swapping},
which is impossible.
\end{proof}

\begin{lemma}
\label{lemma:Hong-Park-08} There is a curve $Z\subset\Pi$ of
degree $k$ such that $|Z\cap\Sigma^{\prime}|\geqslant(2r-1)k+1$.
\end{lemma}

\begin{proof}
Suppose that no $(2r-1)k+1$ points of the set $\Sigma^{\prime}$
lie on a curve of degree $k$ for every natural number $k$. Let us
show that this assumption leads to a contradiction.

The finite subset $\Sigma^{\prime}\setminus P^{\prime}\subset\Pi$
and the natural number $3r-4$ does not satisfy ar least one of the
hypotheses of Theorem~\ref{theorem:Bese}, because every surface in
$\mathbb{P}^{3}$ of degree $3r-4$ containing all points of the set
$\Sigma\setminus P$ must contain the point $P$. However, the
inequalities
$$
\big|\Sigma^{\prime}\setminus P^{\prime}\big|\leqslant
\big(2r-1\big)r\leqslant\mathrm{max}\left(\Big\lfloor\frac{3r-1}{2}\Big\rfloor\Big(3r-1-\Big\lfloor\frac{3r-1}{2}\Big\rfloor\Big),
\Big\lfloor\frac{3r-1}{2}\Big\rfloor^{2}\right)
$$
hold, and at most $3r-4$ points of the set
$\Sigma^{\prime}\setminus P^{\prime}$ can lie on a line, because
$3r-4\geqslant 2r-1$ and at most $2r-1$ points of the set
$\Sigma^{\prime}$  can lie on a line by
Lemma~\ref{lemma:zero-dimensional}.

We see that at least $k(3r-1-k)-1$ points of the set
$\Sigma^{\prime}\setminus P^{\prime}$ lie on a curve of degree~$k$
such that $2\leqslant k\leqslant (3r-1)/2$, which implies that
$k=r$, because at most $k(2r-1)$ points of the set
$\Sigma^{\prime}$ lie on a curve of degree $k$, and
$|\Sigma^{\prime}\setminus P^{\prime}|\leqslant (2r-1)r$. Thus, we
conclude that there is a curve $C\subset\Pi$ of degree $r$ that
contains at least $(2r-1)r-1$ points of $\Sigma^{\prime}\setminus
P^{\prime}$.

The curve $C$ contains $P^{\prime}$, because otherwise there is a
curve in $\Pi$ of degree $3r-4$ that contains the set
$\Sigma^{\prime}\setminus P^{\prime}$ and does not contain the
point $P^{\prime}$. Hence, the curve $C$ contains at least
$(2r-1)r$ points of the set $\Sigma^{\prime}$. Thus, the curve $C$
is reducible by Lemma~\ref{lemma:Hong-Park-07}.

Let $C=\sum_{i=1}^{\alpha}C_{i}$, where $C_{i}$ is an irreducible
curve of degree $d_{i}\geqslant 1$ and $\alpha\geqslant 2$. Then
$$
(2r-1)r\leqslant\big|C\cap\Sigma^{\prime}\big|\leqslant
\sum_{i=1}^{\alpha}\big|C_{i}\cap\Sigma^{\prime}\big|\leqslant\sum_{i=1}^{\alpha}(2r-1)\mathrm{deg}\big(C_{i}\big)=(2r-1)r,
$$
which implies that the curve $C_{i}$ contains $(2r-1)d_{i}$ points
of the set $\Sigma$, and every point of the set $\Sigma$ is
contained in at most one curve $C_{i}$.

Let $C_{\upsilon}$ be the irreducible component of the curve $C$
that contains $P^{\prime}$, and $\Upsilon$ be a subset of the set
$\Sigma$ that contains all points of the set $\Sigma$ that are
mapped to the curve $C_{\upsilon}$ by the projection $\psi$. Then
$|\Upsilon|=(2r-1)d_{\upsilon}$, but less than
$(2r-1)d_{\upsilon}$ points of the set $\Sigma$ lie on a curve of
degree $d_{\upsilon}$. Hence, the points of the set $\Upsilon$
impose independent linear conditions on the homogeneous forms of
degree $3(d_{\upsilon}-1)$ by Lemmas~\ref{lemma:zero-dimensional}
and \ref{lemma:non-vanishing}.

There is a surface $F\subset\mathbb{P}^{3}$ of degree
$3(d_{\upsilon}-1)$ that contains the set $\Upsilon\setminus P$
and does not contain the point $P$. Let $Y_{i}$ be a cone in
$\mathbb{P}^{3}$ over $C_{i}$ with the vertex $O$. Then the
surface
$$
F\cup\bigcup_{i\ne\upsilon} Y_{i}\in\big|\mathcal{O}_{\mathbb{P}^{3}}\big(2d_{\upsilon}-3+r\big)\big|%
$$
contains $\Sigma\setminus P$ and does not contain $P$, but
$2d_{\upsilon}-3+r\leqslant 3r-4$, which is a contradiction.
\end{proof}

There is a disjoint union
$\cup_{j=k}^{l}\cup_{i=1}^{c_{j}}\Lambda_{j}^{i}\subseteq\Sigma$,
where $\Lambda_{j}^{i}$ is a subset of the set $\Sigma$ such that
the inequality $|\Lambda_{j}^{i}|>(2r-1)j$ holds, all points of
the subset $\psi(\Lambda_{j}^{i})$ is contained in an irreducible
curve in $\Pi$ of degree $j$, and at most $(2r-1)t$ points of the
subset
$$
\psi\Big(\Sigma\setminus\Big(\bigcup_{j=k}^{l}\bigcup_{i=1}^{c_{j}}\Lambda_{j}^{i}\Big)\Big)\subsetneq\Sigma^{\prime}\subset\Pi\cong\mathbb{P}^{2}%
$$
lie on a curve in $\Pi$ of degree $t$. Then $k\geqslant 3$ by
Lemma~\ref{lemma:conic}, and $k<r$ by
Lemma~\ref{lemma:Hong-Park-07}.

Put $\Lambda=\cup_{j=k}^{l}\cup_{i=1}^{c_{j}}\Lambda_{j}^{i}$. Let
$\Xi_{j}^{i}$ be the base locus of the linear system of surfaces
in $\mathbb{P}^{3}$ of degree $j$ that contains all points of the
set $\Lambda_{j}^{i}$. Then $\Xi_{j}^{i}$ is a finite set by
Lemma~\ref{lemma:zero-dimensional} and
\begin{equation}
\label{equation:Hong-Park-main-inequality}
\big|\Sigma\setminus\Lambda\big|\leqslant\big(2r-1\big)r+1-\sum_{i=k}^{l}c_{i}\Big(\big(2r-1\big)i+1\Big)\leqslant\big(2r-1\big)\Big(r-\sum_{i=k}^{l}ic_{i}\Big),%
\end{equation}
which implies that $\sum_{i=k}^{l}ic_{i}\leqslant r$.

\begin{remark}
\label{remark:Hong-Park-r-2} The inequality
$\sum_{i=k}^{l}ic_{i}\leqslant r-1$ holds, because the equality
$\sum_{i=k}^{l}ic_{i}=r$ and the
inequalities~\ref{equation:Hong-Park-main-inequality} imply that
$k=l=r$, but $k<r$ by Lemma~\ref{lemma:Hong-Park-07}.
\end{remark}

It follows from Lemma~\ref{lemma:non-vanishing} that the points of
$\Xi_{j}^{i}$ impose independent linear conditions on homogeneous
forms on $\mathbb{P}^{3}$ of degree $3(j-1)$. Put
$\Xi=\bigcup_{j=k}^{l}\bigcup_{i=1}^{c_{j}}\Xi_{j}^{i}$. Then
\begin{equation}
\label{equation:Hong-Park-auxiliary-inequality}
\big|\Sigma\setminus\Big(\Xi\cap\Sigma\Big)\big|\leqslant\big(2r-1\big)r-\sum_{i=k}^{l}c_{i}\big(2r-1\big)i.
\end{equation}

There are surfaces $F$ and $G$ in $\mathbb{P}^{3}$ of degree
$\sum_{i=k}^{l}3(i-1)c_{i}$ and $\sum_{i=k}^{l}ic_{i}$
respectively such that $F$ contains $(\Xi\cap\Sigma)\setminus P$
and does not contain $P$, but $G$ contains $\Xi\cap\Sigma$ and
does not contain any point in $\Sigma\setminus(\Xi\cap\Sigma)$. In
particular, we have $\Sigma\not\subseteq\Xi$, because
$$
\sum_{i=k}^{l}3\big(i-1\big)c_{i}\leqslant\sum_{i=k}^{l}3ic_{i}-3c_{k}\leqslant 3r-6<3r-4.%
$$

Put $\bar{\Sigma}=\psi(\Sigma\setminus(\Xi\cap\Sigma))$ and
$d=3r-4-\sum_{i=k}^{l}ic_{i}$. Then it follows from
Lemma~\ref{lemma:swapping} that there is a point
$\bar{Q}\in\bar{\Sigma}$ such that every curve in $\Pi$ of degree
$d$ that contains $\bar{\Sigma}\setminus\bar{Q}$~must pass through
the point $\bar{Q}$ as well. Therefore, we can not apply
Theorem~\ref{theorem:Bese} to the points of the subset
$\bar{\Sigma}\setminus\bar{Q}\subset\Pi$ and the natural number
$d$.

The proof of Lemma~\ref{lemma:surfaces-IV} implies that the
inequality
$$
\big|\bar{\Sigma}\setminus\bar{Q}\big|\leqslant\big(2r-1\big)\Big(r-\sum_{i=k}^{l}c_{i}i\Big)-1\leqslant\Big\lfloor\frac{d+3}{2}\Big\rfloor^{2}%
$$
holds, but $d=3r-4-\sum_{i=k}^{l}ic_{i}\geqslant 2r-3\geqslant 3$,
because $\sum_{i=k}^{l}ic_{i}\leqslant r-1$, which implies that at
least $t(d+3-t)-1$ points of $\bar{\Sigma}\setminus\bar{Q}$ lie on
a curve in $\Pi$ of degree $t\leqslant (d+3)/2$.

\begin{lemma}
\label{lemma:Hong-Park-09} The inequality $t\ne 1$ holds.
\end{lemma}

\begin{proof}
Suppose that $t=1$. Then at least $d+1$ points of the set
$\bar{\Sigma}\setminus\bar{Q}$ lie on a line, which implies the
inequality $d+1\leqslant 2r-2$ by
Lemmas~\ref{lemma:zero-dimensional} and \ref{lemma:Hong-Park-06}.

The inequality $d+1\leqslant 2r-2$ implies that
$\sum_{i=k}^{l}ic_{i}=r-1$ and $d=2r-3$.

It follows from the
inequality~\ref{equation:Hong-Park-auxiliary-inequality} that
$|\Sigma\setminus(\Xi\cap\Sigma)|\leqslant 2r-1$, which implies
that the points of the set $\Sigma\setminus(\Xi\cap\Sigma)$ impose
independent linear conditions on the homogeneous forms of degree
$2r-3$ by Theorem~2 in \cite{EiJ87}, which is impossible by
Lemma~\ref{lemma:swapping}.
\end{proof}

There is a curve $C\subset\Pi$ of degree $t\geqslant 2$ that
contains at least $t(d+3-t)-1$ points of
the~set~$\bar{\Sigma}\setminus\bar{Q}$, which implies that
$t(d+3-t)-1\leqslant|\bar{\Sigma}\setminus\bar{Q}|$ and
$t(d+3-t)-1\leqslant(2r-1)t$, because at most $(2r-1)t$ points of
the set $\bar{\Sigma}$ lie on a curve of degree $t$.
Therefore,~we~see that $t\geqslant r-\sum_{i=k}^{l}ic_{i}$,
because $t\geqslant 2$. It follows from the
inequalities~\ref{equation:Hong-Park-main-inequality}~that
$$
\big(2r-1\big)\Big(r-\sum_{i=k}^{l}ic_{i}\Big)-1\geqslant\big|\bar{\Sigma}\setminus\bar{Q}\big|\geqslant t(d+3-t)-1\geqslant \Big(r-\sum_{i=k}^{l}ic_{i}\Big)\big(2r-1\big)-1,%
$$
which implies that $t=r-\sum_{i=k}^{l}ic_{i}$, the curve $C$
contains all points of the set $\bar{\Sigma}\setminus\bar{Q}$, and
the inequalities~\ref{equation:Hong-Park-main-inequality} are
actually equalities. Namely, we have $\Sigma\cap\Xi=\Lambda$ and
$$
\big|\Sigma\setminus\Lambda\big|=\big(2r-1\big)r+1-\sum_{i=k}^{l}c_{i}\Big(\big(2r-1\big)i+1\Big)=\big(2r-1\big)\Big(r-\sum_{i=k}^{l}ic_{i}\Big),%
$$
which implies that $l=k$, $c_{k}=1$, $d=3r-4-k$ and
$\sum_{i=k}^{l}ic_{i}=k$.

\begin{lemma}
\label{lemma:Hong-Park-10} The curve $C$ contains all points of
the set $\bar{\Sigma}$.
\end{lemma}

\begin{proof}
Suppose that $C$ does not contain the set $\bar{\Sigma}$. Then $C$
does not contains $\bar{Q}$, which implies that there is a curve
in $\Pi$ of degree $r-k$ that contains the set
$\bar{\Sigma}\setminus\bar{Q}$ does not contain the point
$\bar{Q}$, which is impossible, because $d\geqslant r-k$.
\end{proof}

Thus, the curve $C$ is a curve of degree $r-k$ that contains the
set $\psi(\Sigma\setminus\Lambda)$, which consists of exactly
$(r-k)(2r-1)$ points of the set $\psi(\Sigma)$. On the other hand,
there is an irreducible curve $Z\subset\Pi$ of degree $k$ that
contains all points of the set $\psi(\Lambda)$, which consists of
exactly $k(2r-1)+1$ points of the set $\psi(\Sigma)$. In
particular, we have
$$
\big|\Sigma\big|=\big|\Sigma\setminus\Lambda\big|+\big|\Lambda\big|=\big(r-k\big)\big(2r-1\big)+k\big(2r-1\big)+1=\big(2r-1\big)r+1.
$$

\begin{lemma}
\label{lemma:Hong-Park-11} The curve $C$ is reducible.
\end{lemma}

\begin{proof}
Suppose that $C$ is irreducible. Then the points of the set
$\Sigma\setminus\Lambda$ impose independent linear conditions on
surfaces of degree $3(r-k-1)$ by
Lemmas~\ref{lemma:zero-dimensional},~\ref{lemma:non-vanishing}~and~\ref{lemma:Hong-Park-05},
but the points of the set $\Lambda$ impose independent linear
conditions on surfaces of degree $3(k-1)$ by
Lemmas~\ref{lemma:zero-dimensional}~and~\ref{lemma:non-vanishing},
which implies that the points of the set $\Sigma$ impose
independent linear conditions on homogeneous forms of degree
$3r-4$ by Lemma~\ref{lemma:swapping}.
\end{proof}

Therefore, we have $C=\sum_{i=1}^{\alpha}C_{i}$, where $C_{i}$ is
an irreducible curve of degree $d_{i}$, which implies that
$r-k=\sum_{i=1}^{\alpha}d_{i}$, the curve $C_{i}$ contains
$(2r-1)d_{i}$ points of the set $\bar{\Sigma}$~for~every~$i$, and
every point of the set $\bar{\Sigma}$ is contained in a single
irreducible component of $C$.

\begin{lemma}
\label{lemma:Hong-Park-12} The curve $Z$ contains the point
$P^{\prime}$.
\end{lemma}

\begin{proof}
Suppose that $P^{\prime}\not\in Z$. Let $C_{\upsilon}$ be an
irreducible component of the curve $C$ that contains the point
$P^{\prime}$, and $\Upsilon$ be a subset of the set $\Sigma$ that
contains all points that are mapped to the curve $C_{\upsilon}$ by
the projection $\psi$. Then $\Upsilon$ contains
$(2r-1)d_{\upsilon}$ points.

The points of the set $\Upsilon$ impose independent linear
conditions on the homogeneous forms of degree $3(d_{\upsilon}-1)$
by Lemmas~\ref{lemma:zero-dimensional},~\ref{lemma:non-vanishing}
and \ref{lemma:Hong-Park-05}. Therefore, there is a surface
$F\subset\mathbb{P}^{3}$ of degree $3(d_{\upsilon}-1)$ that
contains $\Upsilon\setminus P$ and does not contain $P$.

Let $Y_{i}$ and $Y$ be the cones in $\mathbb{P}^{3}$ over the
curves $C_{i}$ and $Z$, respectively, whose vertex is the point
$O$. Then the union $F\cup Y\cup\cup_{i\ne\upsilon}Y_{i}$ is a
surface of degree $2d_{\upsilon}-3+r\leqslant 3r-4$ that contains
the set $\Sigma\setminus P$ and does not contain the point $P$,
which is a contradiction.
\end{proof}

The proof of Lemma~\ref{lemma:Hong-Park-12} implies that the
points of the set $\Sigma\setminus\Lambda$ impose independent
linear conditions on homogeneous forms on $\mathbb{P}^{3}$ of
degree $3r-4-k$, but we already know that the points of the set
$\Lambda$ impose independent linear conditions on homogeneous
forms of degree $3(k-1)$ by
Lemmas~\ref{lemma:zero-dimensional}~and~\ref{lemma:non-vanishing},
which is impossible by Lemma~\ref{lemma:swapping}.

\end{document}